\documentclass[10pt]{article}
\usepackage{amsmath}
\usepackage{amssymb}
\usepackage{amsthm}
\usepackage{graphicx}
\usepackage{verbatim}
\usepackage{enumerate}
\font\elevenss=cmss11

\font\eightss=cmss8

\font\sixss=cmss8 at 6pt

\newfam\ssfam
\textfont\ssfam=\elevenss \scriptfont\ssfam=\eightss
 \scriptscriptfont\ssfam=\sixss
\def\ss{\fam\ssfam \elevenss}%

\theoremstyle{plain}
\newtheorem{thm}{Theorem}[section]

\newtheorem{lem}[thm]{Lemma}
\newtheorem{pr}[thm]{Proposition}
\newtheorem{cor}[thm]{Corollary}
\newtheorem{defn}[thm]{Definition}

\theoremstyle{remark}

\newtheorem*{unremark}{Remark}

\def\Z{\mathbb{Z}}

\def\R{\mathbb{R}}

\def\ee{\varepsilon}
\def\E{{\mathbb E}}
\def\P{{\mathbb P}}
\def\Cox{\hfill \Box}

\def\one{{\bf 1}}
\def\|{{\, | \, }}
\def\F{{\mathcal F}}

\def\sgn{\mbox{sgn}\,}
\def\up{\mbox{\ss up}}
\def\down{\mbox{\ss down}}
\def\med{{\hat{\ell}}}

\def\state{{\mathcal S}}

\begin{document}
\begin{titlepage}
\begin{center}
{\large \bf Non-universality for longest increasing subsequence of
a random walk}
\\[5ex]
\end{center}

\begin{flushright}
Robin Pemantle\footnote{Department of Mathematics,
University of Pennsylvania,
209 South 33rd Street,
Philadelphia, PA 19104, USA, pemantle@math.upenn.edu}$^,$
\footnote{Research supported in part by NSF grant \# DMS-1209117},
and
Yuval Peres\footnote{Microsoft Research,
1 Microsoft Way, Redmond, WA, 98052, USA,
peres@microsoft.com},
\end{flushright}

\vfill

\noindent{\bf Abstract:}
The longest increasing subsequence of a random walk with mean zero
and finite variance is known to be $n^{1/2 + o(1)}$.
We show that this is not universal for symmetric random walks.
In particular, the symmetric Ultra-fat tailed random walk has a longest
increasing subsequence that is asymptotically at least $n^{0.690}$
and at most $n^{0.815}$.  An exponent strictly greater than $1/2$
is also shown for the symmetric stable-$\alpha$ distribution
when $\alpha$ is sufficiently small.
\vfill

\\[2ex]


\end{titlepage}

\section{Introduction}
\label{sec:intro}

It is well known that the longest increasing subsequence (LIS) of
a sequence of $n$ IID non-atomic random variables has length $(2 + o(1)) \sqrt{n}$ with high probability.
(see~\cite{vershik-kerov,logan-shepp}).
A different model was considered by Angel et al.~\cite{angel-balka-peres}.
Let $S_n := \sum_{k=1}^n X_k$ be the partial sums of random walk with
mean zero and finite variance.  Angel et al.\ show that the LIS of the
partial sum sequence $(S_1 , \ldots , S_n)$ has length $n^{1/2 + o(1)}$.
They do not shed any light on what happens when the finite variance
hypothesis is removed.  When the second moment, and possibly the first,
are undefined, it makes sense to consider other ways to keep the
walk from having a drift.  Here we consider random walk trajectories
whose increments are symmetric about zero.  We show that
random walks whose increments have fat tails will have a longer LIS
than do those with finite variance.

The cleanest model in which this occurs is the so-called Ultra-fat tailed
distribution, which is a distribution not on $\R$ but on a
non-archimedean totally ordered space $\state$ described in
Section~\ref{sec:fat} below.  There, we are able to show that
the LIS has length at least $n^{0.69}$; see Theorem~\ref{th:main}
below.  However, the result also holds for real random walks with
fat tails, such as the symmetric stable-$\alpha$ when $\alpha$
is sufficiently small.  We also show that the LIS has length
at most $n^{0.82}$.  Neither of these exponents is believed to
be sharp, however empirical studies suggest that the LIS exponent
for the Ultra-fat tailed distribution is roughly $0.72$, so nearer
to our lower bound.  Numerical evidence also suggests that for
stable laws, the exponent varies, interpolating between this and $1/2$.

The organization of the rest of the paper is as follows.
The next section contains definitions, notation and preliminary
facts.  Section~\ref{sec:lower} proves the $n^{\beta_0 + o(1)}$ lower bound
with an explicit constant $\beta_0$ slightly larger than $0.69$.
Section~\ref{sec:upper} proves the $n^{\beta_1 + o(1)}$ upper
bound, with an explicit constant $\beta_1$ slightly less than $0.9$.
Section~\ref{sec:real} extends the lower bound from the Ultra-fat tail case
to actual fat-tailed distributions.  We conclude with some further
remarks and questions.

\section{Definitions and results}
\label{sec:fat}

\subsection{Ultra-fat tailed distribution}

We begin by defining the state space $\state$, which is a free
$\Z$-module with one generator $x$ for each $x \in (0,1)$.
In other words, elements of $\state$ are finite formal
linear combinations of the symbols $\{ x : 0 < x < 1 \}$
with coefficients in $\Z$.  There should be no confusion
between the formal symbol $x$ and the real number $x$
as coefficients take only integer values and are
always written on the left.

Endow $\state$ with the lexicographic order.  Formally, if
$\alpha = \sum_{x \in F} a_x x$ and $\beta = \sum_{x \in G} b_x x$,
we may define this order relation by induction on the minimum length
$m \wedge n$ of $\alpha$ and $\beta$ as follows.  For
$\alpha  = \sum_{x \in F} a_x x \in \state$, define its degree
by $|\alpha| := \sup \{ t : a_t \neq 0 \}$.  By convention $|0| = 0$.
We define comparisons to~0 by $\alpha > 0$ if and only if
$|F| > 0$ and $a_{|\alpha|} > 0$ and $\alpha < 0$ if and only if
$-\alpha > 0$.  For elements $\alpha = \sum_{x \in F} a_x x $
and $\beta = \sum_{x \in G} b_x x$, assuming $|F| , |G| > 0$,
inductively define $\alpha > \beta$ if and only if one of
the following conditions holds.
\begin{enumerate}[(i)]
\item $|\alpha| = t > |\beta|$ and $a_t > 0$;
\item $|\beta| = t > |\alpha|$ and $b_t < 0$;
\item $|\alpha| = t = |\beta|$ and $a_t > b_t$;
\item $|\alpha| = t = |\beta|$ and $a_t = b_t$ and
   $\alpha - a_t t > \beta - b_t t$.
\end{enumerate}
This defines a total order on $\state$ consistent with addition:
$\alpha > \beta$ and $\gamma \geq \delta$ implies $\alpha + \gamma
\geq \beta + \delta$.

\subsection{Ultra-fat tailed random walk}

Define $F : [-1,1] \to \state$ by $F(0) = 0$ and $F(x) = \sgn (x) |x|$
for $x \neq 0$.  Let $\{ U_n : n \geq 1 \}$ be an IID collection
of real random variables uniform on $[-1,1]$.  The law of $F(U_1)$
is called the {\bf Ultra-fat tailed distribution}.  Let $X_n := F(U_n)$
and $S_n := \sum_{k=1}^n X_k$.  The sequence $\{ S_n : n \geq 0 \}$
is called the so-called {\bf Ultra-fat tailed random walk}.  A sequence
$n_1 < \cdots < n_k$ is an increasing subsequence if
$S_{n_i} < S_{n_j}$ for all $1 \leq i < j \leq k$.  The Ultra-fat tailed
distribution has been used elsewhere, without a formal definition;
see, e.g.,~\cite{limic-pemantle}.  The following theorem is the
main result of this note.

\subsection{Main result}

Let $\{ S_n \}$ be a random walk on $\state$ with increments
from the Ultra-fat tailed distribution.  Let $L(t)$ denote the
length of the LIS of $(S_1 , \ldots , S_t)$.  We remark that,
by convention, we have not allowed $S_0$ to be an element of
the LIS, hence the increment $X_1$ will play no role.

\begin{thm} \label{th:main}
There are numbers $1/2 < \gamma < \delta < 1$  such
that as $t \to \infty$,
$$\P (  t^\gamma \ge L(t) \le t^\delta) \to 1 \, .$$
In particular, one can take $\gamma = 0.690$ and $\delta = 0.815$.
\end{thm}

\begin{unremark}
It can be shown (see Section~\ref{sec:questions}) that neither
exponent is sharp.
\end{unremark}

\section{Preliminary results}

We look at the growth rate of various deterministic random
functions going between the time variable for the random walk
and the length variable for the corresponding LIS.  Because
of the   proliferation of notation, we will organize
by using $t$ and nearby letters when possible for quantities
in the time domain and $\ell$ and nearby letters for quantities
in the length domain.  As usual, we use upper case letters such
as $L$ and $T$ for random quantities.

The random variables $\{ L(t) : t \geq 1 \}$ have already
been introduced and follow this notational scheme.
For $\ell \geq 1$ let $T(\ell) := \inf \{ t : L(t) \geq \ell \}$
denote the random time that the LIS first reaches length $\ell$.
Thus $L(T(\ell)) = \ell$ and $T(L(t)) \leq t$.

The magnitudes of the steps are the values $|U_1| , |U_2|, \ldots$.
Trivially, the order type of the first $t$ of these is uniform
on all $t!$ possible orders and independent of the sign vector,
which is also uniform on $\{ \pm 1 \}^t$.  This allows for the
usual conditioning identities.  For example, if the variable of
the greatest magnitude is $U_\sigma$ then the order types of
$(|U_1| , \ldots , |U_{\sigma - 1}|)$ and $(|U_{\sigma + 1}| ,
\ldots , |U_t|)$ are independent and uniform.
Also immediate is the following Markov property.  Construct
the random variables $\{ U_n \}$ as the coordinate functions
on the canonical space $\Omega := [-1,1]^\infty$ with normalized
Lebesgue measure.  Let $\theta : \Omega \to \Omega$ be the
shift $(U_1, U_2, U_3, \ldots) \mapsto (U_2, U_3, \ldots)$.
Let $\F_t := \sigma (U_1 , \ldots , U_t)$ and let $\tau$
be a stopping time with respect to the filtration $\{ \F_t \}$.
Then conditionally on $\F_\tau$, the sequence $\{ X_{\tau + n} \}$
is distributed as the unconditional sequence $\{ X_n \}$.

This is all pretty trivial but it allows us to state two
important relationships, one sub-additive and super-additive:
\begin{eqnarray}
L(s+t) & \leq & L(s) + L(t) \circ \theta^s \label{eq:sub}
   \, ; \\[1ex]
T(\ell + m) & \geq & T(\ell) + T(m)  \circ \theta^{T(\ell)}
   \, . \label{eq:super}
\end{eqnarray}
Intuitively, the first of these holds because any increasing subsequence
of $(S_1, \ldots , S_{s+t})$ has at most $L(s)$ entries in $[s]$
and $L(t) \circ \theta^s$ elements in $\{ s+1 , \ldots , s+t \}$.
The second holds because to get an increasing subsequence of
length $\ell + m$ one first needs one of length $\ell$, and must
then find one of length $m$ among the remainder of the sequence.
These properties do not rely on the Ultra-fat tailed distribution and
hold for the LIS of any random walk.

\begin{defn}[NBU]
Say that a random variable $X$ is new better than used (NBU) if
for every pair of positive integers $a$ and $b$,
$$\P (X \geq a+b) \leq \P(X \geq a) \P (X \geq b) \, .$$
\end{defn}
\noindent{The} terminology comes from reliability
theory~\cite{BMP63,barlow-proschan}, where the
inequality rewritten as $\P (X \geq a+b | X \geq a) \leq \P (X \geq b)$
says that a new light bulb has a better chance of surviving $b$ units
of time, than does a light bulb that has been used for $a$ units of time.

We recall a basic property of NBU variables.
\begin{lem} \label{lem:quantile}
If the random variable $X$ is NBU and $\P(X<q)=\epsilon$, then $\E[X] \le  q/\epsilon$.
\end{lem}
\noindent{\sc Proof:}  The NBU assumption implies that $\P(X \ge kq) \le (1-\epsilon)^k$, so   $X/q$ is stochastically dominated by a Geometric Variable of mean $1/\epsilon$.

\begin{pr} \label{pr:NBU}
For LIS of any random walk, each random variable $L(t)$ is NBU.
\end{pr}

\noindent{\sc Proof:} This follows from~\eqref{eq:sub}
and~\eqref{eq:super}.  The event $\{ L(t) \geq a+b \}$ is
the intersection of the events $\{ T(a) \leq t \}$ and
$\{ T(b) \circ \theta^{T(a)} \leq t - T(a) \}$.  Because
$T(b) \circ \theta^{T(a)}$ is independent of $\F_{T(a)}$,
\begin{eqnarray*}
\P \bigl(L(t) \geq a+b \bigr) & = & \P (T(a) \leq t)  \, \cdot \,
    \P \Bigl( T(b) \circ \theta^{T(a)} \leq t - T(a)  \Bigr)    \\[1ex]
& \le & \P (L(t) \geq a) \, \cdot \, \P (L(t) \geq b) \, .
\end{eqnarray*}
$\Cox$

\begin{pr}
Let $X$ be NBU with mean $\mu$ and let $Y$ be geometric started from zero
with mean $\mu$ (that is, one less than a geometric of mean $\mu + 1$).
Denote
\begin{eqnarray*}
a_n & := & \P (X \geq n) \\
A_n & := & \sum_{k=n}^\infty a_k \\
g_n & := & \P (Y \geq n)
   = \left ( \frac{\mu}{1+\mu} \right )^n \\
G_n & := & \sum_{k=n}^\infty g_k = G_0 g_n \\
\end{eqnarray*}
where $G_0 = 1 + \mu$.  Then for all $n$ we have $A_n \leq G_n$.
\end{pr}

\noindent{\sc Proof:}  Let $t$ be the least integer such that
$a_t < g_t$.  Then $t$ is at least~1 because $a_0 = g_0 = 1$.
Also $t$ is finite unless $X$ and $Y$ have the same distribution
because
$$\sum_{n=0}^\infty a_n = \sum_{n=0}^\infty g_n = 1 + \mu \, .$$
Suppose first that $n \leq t$.  Then
$$A_n = 1 + \mu - \sum_{k=0}^{n-1} a_k \leq
   1 + \mu - \sum_{k=0}^{n-1} g_k = G_n \, .$$
Now suppose that $n > t$ and assume for induction that
$A_m \leq G_m$ for all $m < n$.  Then using the NBU property
and induction,
$$A_n = \sum_{k \geq n} a_k \leq a_t \sum_{k \geq n-t} a_k
   = a_t A_{n-t} \leq g_t G_{n-t} = G_n \, ,$$
completing the induction.
$\Cox$

\begin{cor} \label{cor:convex domination}
Let $X$ be NBU and $Y$ be geometric started from zero with the same mean.
Then for any convex function $\phi$,
$$\E \phi (X) \leq \E \phi (Y) \, .$$
\end{cor}

Remark: An equivalent conclusion is that $X$ has the distribution of
some conditional expectation of $Y$.

\noindent{\sc Proof:} As before, let $a_n$ and $A_n$ be tail probabilities
for $X$ and their tail sums.  Letting $\Delta h (n)$ denote
$h (n+1) - h (n)$, we sum by parts twice to obtain
$$\E \phi (X) = \phi (0) + \sum_{n=1}^\infty \Delta \phi (n-1) a_n
   = \phi (0) + \Delta \phi (0) A_1 + \sum_{n=2}^\infty \Delta \Delta
   \phi (n-2) A_n \, .$$
Similarly,
$$\E \phi (Y) = \phi (0) + \Delta \phi (0) G_1 + \sum_{n=1}^\infty
   \Delta \Delta \phi (n-2) G_n \, .$$
Because $\phi$ is convex, $\Delta \Delta \phi (k) \geq 0$ for $k \geq 0$.
Together with $A_k \leq G_k$ for all $k$ and $A_1 = G_1 = \mu$, this
proves the corollary.
$\Cox$

\section{Proof of the lower bound}
\label{sec:lower}

In this section we prove the lower bound in Theorem~\ref{th:main}.
In terms of universality, this direction is the more interesting,
as it shows the Ultra-fat tailed walk to be in a different LIS-universality
class from mean zero finite variance walks.

For $1 \leq m < n$, let $L(m,n) := L(n-m) \circ \theta^m$,
in other words, it is the length of the LIS of $(S_m , \ldots ,
S_n)$ (recall that, by convention, the LIS cannot include the
initial element, $S_m$).  Of course
$L(m,n)$ has the same distribution as $L(n-m)$.  Define
$\sigma(n)$ to be the almost surely unique $k \in [n] \setminus \{ 1 \}$
such that $|U_k| = \max_{2 \leq j \leq n} |U_j|$.  In other words,
$\sigma(n)$ is the time at which the random walk completed
its largest magnitude step among those occuring after time~1
and before time~$n$.  Let $\up_n := \{ U_{\sigma(n)} > 0 \}$
denote the event that this greatest magnitude increment was positive.
The complementary event is denoted $\down_n$.  On $\up_n$, one
has the inequality $S_j > S_i$ whenever $j \geq \sigma(n) > i$.
Therefore, the increasing subsequences of $[n]$ are precisely
the unions $A \cup B$ where $A$ is an increasing subsequence
of $[\sigma(n) - 1]$ and $B$ is an increasing subsequence of
$[n] \setminus [\sigma(n) - 1]$.  It follows that
$L(n) = L(\sigma(n)-1) + L(\sigma (n)-1,n)$
On $\down_n$, one has $S_j < S_i$ whenever $j \geq \sigma(n) > i$,
hence the increasing subsequences of $[n]$ are precisely the sets
that are either an increasing subsequence of $[\sigma(n) - 1]$ or
of $[n] \setminus [\sigma(n) - 1]$.  We have therefore proved:
\begin{pr} \label{pr:pointwise}
The sequence of random variables $\{ L(n) \}$ satisfies the recursion
\begin{eqnarray*}
L(n) & = & \;\;\;\;
   \one_{\up_n} \left [ L(\sigma(n) - 1) + L(\sigma(n) - 1,n) \right ] \\
&& + \; \one_{\down_n}
   \max \left \{ L(\sigma(n) - 1) , L(\sigma(n) - 1,n) \right \} \, .
\end{eqnarray*}
$\Cox$
\end{pr}

We will prove the lower bound in Theorem~\ref{th:main}
by using the recursion to obtain the following lower bound
on $\E L(n)$.
\begin{lem} \label{lem:E LB}
Let $\beta_0$ be the positive solution to $x + 2^{-1-x} = 1$,
whose decimal expansion begins $0.690069$.  Then
$\E L(n) \geq n^{\beta_0 - o(1)}$.
\end{lem}

\noindent{\sc Proof:}
Let $a_n = \E L(n)$.  Note that, conditional on $\sigma(n)$, the three 
random variables $L(\sigma(n) - 1)$, $L(\sigma(n)-1,n)$ and
$\one_{\up_n}$ are all independent and distributed respectively as
$L(k-1), L(n-k+1)$ and $\mbox{Bernoulli}(1/2)$, where
$k = \sigma(n)$.  Therefore, using $L(\max\{a,b\})$ as a lower
bound for $\max \{ L(a) , L(b) \}$ in the second line, we have
\begin{eqnarray}
a_n & = & \frac{1}{2(n-1)} \sum_{k=2}^n \left ( a_{k-1} + a_{n-k+1} \right )
   + \frac{1}{2(n-1)} \sum_{k=2}^n \E \max \{ L(k-1) , L(n-k+1 \}
   \nonumber \\
& \geq & \frac{1}{n-1} \sum_{k=1}^{n-1} a_k
   + \frac{1}{n-1} \sum_{k=n/2}^{n-1}
   a_k  \left ( 1 - \frac{1}{2} \delta_{k,n/2} \right ) \, .
   \label{eq:recursion}
\end{eqnarray}

The key observation is that for $\beta < \beta_0$ and
sufficiently large $n$,
\begin{equation} \label{eq:reverse}
n^\beta \leq \frac{1}{n-1} \sum_{k=1}^{n-1} k^\beta
   + \frac{1}{n-1} \sum_{k=n/2}^{n-1}
   k^\beta  \left ( 1 - \frac{1}{2} \delta_{k,n/2} \right ) \, .
\end{equation}
Indeed, dividing~\eqref{eq:reverse} through by $n^\beta$, the
right-hand side is a Riemann sum approximation for
$$c_\beta := \int_0^1 x^\beta \, dx + \int_{1/2}^1 x^\beta \, dx$$
which evaluates to
$$\frac{1}{\beta + 1} \left ( 2 - 2^{-\beta - 1} \right ) \, .$$
As a function of $\beta$, the quantity $c_\beta$ decreases as $\beta$
varies over $[0,1]$, passing through the value~1 at $\beta = \beta_0$.
Therefore, for $\beta < \beta_0$, we have
$$\frac{1}{n-1} \sum_{k=1}^{n-1} k^\beta
   + \frac{2}{n-1} \sum_{k=n/2}^{n-1}
  k^\beta \left ( 1 - \frac{1}{2} \delta_{k,n/2} \right )
= n^\beta \left [ c_\beta - o(1) \right ] > n^\beta$$
provided that $n > N(\beta)$, where $N(\beta)$ is sufficiently large
so that the $o(1)$ term is less than $c_\beta - 1$.

The rest is easy.  Fixing $\beta < \beta_0$, we may pick $C = C(\beta)$
such that $a_n \geq C n^\beta$ for all $n \leq N(\beta)$.  We claim, by
induction, that this is true for all $n > N(\beta)$ as well.
Indeed, assuming it to be true for $n-1$, we see that the
right-hand side of~\eqref{eq:recursion}, which is a lower bound
for $a_n$, is at least $C(\beta)$ times the right-hand side
of~\eqref{eq:reverse}.  Because $n > N(\beta)$, we see
from~\eqref{eq:reverse} that this is at least $C n^\beta$,
proving the claim.

We have shown that for all $\beta < \beta_0$ there exists a $C$
such that  $\E L(n) \geq C n^\beta$ for all $n$.  This completes
the proof of Lemma~\ref{lem:E LB}. 
$\Cox$

\noindent{\sc Proof of lower bound in Theorem~\protect{\ref{th:main}}:}
Fix $\gamma<\beta<\beta_0$. The preceding lemma gives $\E L(n) \geq C n^\beta$.    By Lemma~\ref{lem:quantile},  we have
$\P[L(n)<n^\gamma] \le n^{\gamma-\beta} \to 0$.
$\Cox$

\section{Proof of upper bound}
\label{sec:upper}

The proof of the upper bound in Theorem~\ref{th:main} is analogous
to the proof Lemma~\ref{lem:E LB} but in the reverse direction.
It reduces to the following result.
\begin{lem} \label{lem:E UB}
Let $\beta_1$ be the positive solution to
$$\frac{2}{1+\beta} - \int_0^{1/2} \frac{x^\beta (1-x)^\beta}
   {x^\beta + (1-x)^\beta} = 1 \, ,$$
whose decimal expansion begins $0.814834$.  Then $\E L_n \leq
n^{\beta_1 + o(1)}$.
\end{lem}

Before proving this, we record the following lemma.
\begin{lem} \label{lem:min}
Let $X_1$ and $X_2$ be independent random variables both of which
are NBU.  Let $a := \E X_1 \leq b := \E X_2$.  Then
$$\E (X_1 \wedge X_2) \geq \frac{ab}{a+b+1} \, .$$
\end{lem}

\noindent{\sc Proof:} Let $Y_1$ and $Y_2$ be independent geometric random variables
(taking values $0,1,2,\ldots$)  with respective means $a$ and $b$.
Corollary~\ref{cor:convex domination} applied to $-(X_1 \wedge s)$
shows that $\E (X_1 \wedge s) \geq \E (Y_1 \wedge s)$ for each fixed $s$.
It follows that
$$\E (X_1 \wedge X_2 \| X_2) \geq \E (Y_1 \wedge X_2 \| X_2)$$
and hence that $\E (X_1 \wedge X_2) \geq \E (Y_1 \wedge X_2)$.
Similar reasoning shows that $\E (Y_1 \wedge X_2) \geq \E (Y_1 \wedge Y_2)$.
This last quantity may be computed exactly.  This is one less than
the minimum of two geometrics (started from~1) with respective success
probabilities $1/(a+1)$ and $1/(b+1)$, which means a combined
success probability of $(a+b+1)/(ab+a+b+1)$.  One less than
the mean is $ab/(a+b+1)$, proving the lemma.
$\Cox$

\noindent{\sc Proof of Lemma~\protect{\ref{lem:E UB}}:}
Again let $a_n$ denote $\E L(n)$.
Fix $n$ and again let $\sigma = \sigma (n)$ denote the time
of the largest magnitude step up to time $n$.
The identity $\max \{ a , b \} = a + b - \min \{ a , b \}$ gives
\begin{equation} \label{eq:L rec}
L(n) = L(\sigma - 1) + L(\sigma-1,n )
   - \one_{\down_n} \, \min \{ L(\sigma - 1) , L(\sigma - 1,n) \} \, .
\end{equation}
The random variables $\one_{\down_n}$, $L(\sigma - 1)$ and
$L(\sigma - 1,n)$ are conditionally independent given $\sigma$.
Now use Lemma~\ref{lem:min} with $X = L(k)$ and $Y = L(n-k)$
where $k$ is the minimum of $\sigma - 1$ and $n - \sigma + 1$.
This gives
$$\E \left [ \one_{\down_n}
   \min \{ L(\sigma) - 1 , L(\sigma - 1,n) \} \right ] \geq
   \frac{1-o(1)}{2} \;
   \frac{\E L(k) \cdot \E L(n-k)}{\E [L(k) + L(n-k)]}$$
where the $o(1)$ term is uniform in $k$ as $n \to \infty$,
coming from the ratio of $a+b+1$ and $a+b$ when $a = \E L(k)$
and $b = \E L(n-k)$.
Plugging this in to~\eqref{eq:L rec} after taking expectations gives
\begin{equation} \label{eq:opp}
a_n \leq \frac{2}{n-1} \sum_{k=2}^n a_k -
   \frac{2}{n-1} \sum_{k=2}^{n/2} \frac{1-o(1)}{2}
   \frac{\E L(k) \cdot \E L(n-k)}{\E [L(k) + L(n-k)]}  \, .
\end{equation}

Again we play the trick of replacing $a_k$ by $k^\beta$ and
approximating the sum by an integral.  Pulling out a factor of
$n^\beta$, the right hand side becomes
\begin{equation} \label{eq:recur}
n^\beta \left [ 2 \int_0^1 x^\beta \, dx -
   \int_0^{1/2} \frac{x^\beta (1-x)^\beta}
    {x^\beta + (1-x)^\beta} \right ) \, dx + o(1) \, .
\end{equation}
The expression~\eqref{eq:recur} is decreasing on $[0,1]$ and passes
through the value~1 at $\beta_1$.  Now fix $\beta > \beta_1$, let
$N(\beta)$ be large enough so that~\eqref{eq:recur} is less than $n^\beta$
for all $n \geq N(\beta)$.  Choosing $C$ so that $a_n \leq C n^\beta$
for $n \leq N(\beta)$, the integral approximation then shows by
induction that $a_n \leq C n^\beta$ for all $n$, finishing the
proof of Lemma~\ref{lem:E UB}. Invoking Markov's inequality then ayields the upper bound in
Theorem~\ref{th:main}.
$\Cox$

\section{Real random walks}
\label{sec:real}

For whose who don't accept the Ultra-fat tailed distribution as a true
random walk, we include the following result.

\begin{thm} \label{th:real}
Let $\{ S_n \}$ be the partial sums of a symmetric stable walk
with index $\alpha$.  For any $\gamma < \beta_0$ there are real
$\alpha , C > 0$ such that the length $L(n)$ of the LIS of the
symmetric stable walk to time $n$ has expectation at least
$C n^\gamma$.
\end{thm}

We begin with a lemma.

\begin{lem} \label{lem:UI}
Fix any index $\alpha \in (0,2)$.  Let $\med (n) = \med_\alpha (n)$
denote the median of $L(n)$ for the symmetric stable walk of index
$\alpha$.  Then the family $\{ L(n) / \med (n) : n \geq 1 \}$ is
uniformly integrable.  In particular,
$$\E \left [ \frac{L(n)}{\med (n)} \; \one_A \right ] \leq g(\P(A))$$
for some function $g$ with $\lim_{\ee \downarrow 0} g(\ee) = 0$.
\end{lem}

\noindent{\sc Proof:} Subadditivity~\eqref{eq:sub} holds
for any random walk.  Therefore, $L(n) / \med (n)$ has uniformly
exponential tails, and uniform integrability follows.
$\Cox$

\noindent{\sc Proof of Theorem}~\ref{th:real}:
Let $W_n := \max \{ |X_k| : 1 \leq k \leq n \}$ and
$Z_n := \sum_{k=1}^n |X_k|$.  Let $\up_n'$ be the event
that $\up_n$ occurs and $W_n > Z_n - W_n$.  On $\up_n'$,
the recursion in Proposition~\ref{pr:pointwise} is
satisfied at $n$.  We will show that
\begin{equation} \label{eq:dwarf}
\P (\up_n \setminus \up_n') \leq g(\alpha)
\end{equation}
for some function $g$ such that $\lim_{t \downarrow 0} g(t) = 0$,
uniformly in $n$.  Assuming this, we can complete the analysis
by showing the alteration to the conclusion of
Proposition~\ref{pr:pointwise} is sufficiently small.

On $\down_n$, the inequality is favorable:
we can still choose to use only the longer segment, hence
$L(n) \geq L((\sigma(n) - 1) \vee (n - \sigma(n)))$.  On $\up_n$
the inequality goes the wrong way, but the difference
is bounded above by $(L(\sigma(n) - 1) + L(\sigma(n)-1,n)) (\one_{\up_n}
- \one_{\up_n'})$.  Assuming~\eqref{eq:dwarf}, we take
expectations, yielding
$$\E L_n \geq \frac{1}{n-1}\sum_{k=1}^{n-1} \E L_k + \frac{1}{n-1} \sum_{k=n/2+1}^{n-1} \E L_k
   - \E \Bigl(L(\sigma(n) - 1) + L(\sigma(n)-1,n)\Bigr) (\one_{\up_n} - \one_{\up_n'}) \, .$$
The subtracted term $\E \Bigl(L(\sigma(n) - 1) + L(\sigma(n)-1,n)\Bigr) (\one_{\up_n}
- \one_{\up_n'})$ is $o(\med (n))$ by~\eqref{eq:dwarf} and
Lemma~\ref{lem:UI}.  The approximation by Riemann sum and
the resulting inequality then finish the proof as before.

It remains to show~\eqref{eq:dwarf}.  We remark that this is the
only place we use specific properties of the distribution other than
symmetry.  The conclusion of the theorem will therefore hold for
any symmetric distribution satisfying~\eqref{eq:dwarf}.  In particular,
there are many more extreme distributions, such as $Z = R e^X$
where $X$ is Cauchy and $R$ is Rademacher, for which the ratio
of the greatest of $n$ picks to the sum of the magnitudes of the
other $n-1$ goes to infinity in probability.  For such distributions,
the conclusion holds for all $\gamma < \beta_0$.

Recall that the symmetric stable variable $X_\alpha$ may be
constructed as the difference of IID positive stable variables
$Y - Z$, each of which is the sum of the points of a Poisson
process with intensity $x^{-1-\alpha}$ on $\R^+$.  For $t \geq 1$,
let $Y_t$ and $Z_t$ denote the sum of the $t^{th}$ power of these
points.  Then $Y_t - Z_t$ is a symmetric stable of index $\alpha/t$.
This coupling of $X_{\alpha'}$ for all $\alpha' \leq \alpha$,
together with the fact that the magnitudes of the Poisson points
are almost surely summable and distinct, shows that $Z_n / W_n
\to 1$ almost surely for fixed $n$ as $\alpha \downarrow 0$.  This
of course implies convergence in probability, so the only thing
remaining to check is uniformity in $n$.

This follows from tightness of two families: $\{ W_n / Z_n \}$ and
$\{ Z_n / (Z_n - Z_n') \}$ where $Z_n'$ is the second highest magnitude
of a Poisson summand.  These both follow elementarily from properties
of the Poisson process of intensity $2 x^{-\alpha-1} \, dx$ on $\R^+$.
Set $\alpha = 1$.  Given $\ee > 0$, choose $K$ such that both
$\P (W_n / Z_n > K)$ and $\P (Z_n / (Z_n - Z_n') > K)$ are
less than $\ee/2$.  When both inequalities are satisfied we have
$W_t / Z_t < 1 + (K-1) (1 - 1/K)^t$.  Choosing $t := t(\ee)$
large enough to make this less than~2, we see that $\alpha \geq t$
makes $\P (\up_n \setminus \up_n') \leq \ee$.  This
proves~\eqref{eq:dwarf} with $g$ the inverse function to $t(\ee)$,
completing the proof of Theorem~\ref{th:real}.
$\Cox$

\section{Further remarks and questions}
\label{sec:questions}

One natural question is to prove that the exponent
$\lim \log L(n) / \log n$ exists.

Another is whether we can obtain better bounds on the exponent
by finding a functional form for the distribution which yield
an inequality when passed through the recursion.

Neither exponent $\beta_0$ nor $\beta_1$ is sharp.
The proof of
$$\liminf \frac{\log \E L(n)}{\log n} \geq \beta_0$$
in fact computes the correct exponent, namely $\beta_0$,
for the length of the greedy increasing subsequence.
The GIS is defined by splitting the sequence at the
location $\sigma$ of the maximum step, and, if the step
is a downward step, throwing away the smaller interval rather
than the interval with the shorter LIS.  The length $Z_n$ of
the greedy increasing subsequence obeys the recursion of
Proposition~\ref{pr:pointwise} but with the max taken on the
inside.  This seems likely to give an exponent not too far from
the correct exponent, but it gives up a non-negligible amount
in the recursion and cannot be sharp.

The proof of
$$\limsup \frac{\log \E L(n)}{\log n} \leq \beta_1$$
does not, as far as we know compute anything natural.
This bound could be improved by finding the correct
function $\phi (n,k)$ that computes a better lower
bound on $\E \min \{ L(k) , L(n-k) \}$.  Lemma~\ref{lem:min}
is best possible assuming only the NBU property, as
the geometric random variable is the extreme case.
However, we know more about $L(k)$.  For example, when
$j$ and $k/j$ are integers then
$$\P (L(k) < \ee \E L(k)) \leq \P (L(k/j) < \ee \E L(k))^j \, .$$
If $\E L(j) = j^{\beta + o(1)}$ then taking choosing $j$ so that
$\ee = c j^{-\beta}$ makes $\P (L(k/j) < \ee \E L(k)) < 1/2$
and results in
$$\P (L(k) < \ee \E L(k)) \leq 2^{-\ee^{(1+o(1))/\beta}} \, .$$
The lower tails on $L(k)$ are thus expected to be very small;
this ought to lead to a better lower bound on $\E L(k) \wedge L(n-k)$,
hence a better exponent in Lemma~\ref{lem:E UB} and in Theorem~\ref{th:main}.

\bibliographystyle{alpha}
\bibliography{RP}

\end{document}